\documentclass[12pt]{amsart}
\usepackage{amscd,amssymb}
\usepackage[graph,frame,poly,arc]{xy}  
\usepackage[plainpages,backref,urlcolor=blue]{hyperref}

\theoremstyle{plain}
\newtheorem{thm}[subsection]{Theorem}
\newtheorem{lem}[subsection]{Lemma}
\newtheorem{prop}[subsection]{Proposition}
\newtheorem{cor}[subsection]{Corollary}

\theoremstyle{definition}
\newtheorem{rk}[subsection]{Remark}
\newtheorem{definition}[subsection]{Definition}
\newtheorem{ex}[subsection]{Example}
\newtheorem{conj}[subsection]{Conjecture}

\numberwithin{equation}{section}
\newcommand{\OO}{{\mathcal O}}

\newcommand{\C}{\mathbb{C}}
\newcommand{\PP}{\mathbb{P}}

\begin{document}

\title [Free divisors and rational cuspidal plane curves]
{Free divisors and rational cuspidal plane curves}

\author[Alexandru Dimca]{Alexandru Dimca$^1$}
\address{Univ. Nice Sophia Antipolis, CNRS,  LJAD, UMR 7351, 06100 Nice, France. }
\email{dimca@unice.fr}

\author[Gabriel Sticlaru]{Gabriel Sticlaru}
\address{Faculty of Mathematics and Informatics,
Ovidius University,
Bd. Mamaia 124, 900527 Constanta,
Romania}
\email{gabrielsticlaru@yahoo.com }
\thanks{$^1$ Partially supported by Institut Universitaire de France.}



\begin{abstract} A characterization of freeness for plane curves in terms of the Hilbert function of the associated Milnor algebra is given as well as many new examples of rational cuspidal curves which are free. Some stronger properties are stated as conjectures.

\end{abstract}
 
\maketitle


\section{Introduction} \label{sec:intro}

A (reduced) curve $C$ in the complex projective plane $\PP^2$ is called { \it free}, or a free divisor,
if the rank two vector bundle $T\langle C\rangle=Der(-logC)$ of logarithmic vector fields  along $C$ splits as a direct sum of two line bundles on $\PP^2$. Note that $C$ is a free divisor exactly when the surface  given by the cone over $C$ is free at its vertex, the origin of $\C^3$, in the sense of K. Saito who introduced this fundamental notion in \cite{KS}. See also \cite{BEG}.

A curve $C$ in  $\PP^2$ is called { \it cuspidal } if all of its singular points $p$ are generalized cusps, i.e. the analytic germs $(C,p)$ are irreducible. If $d$ is the degree of $C$ and $m$ is the maximal multiplicity of the singular points of $C$, then $C$ is called a curve of type $(d,m)$. Let $\kappa$ be the total number of cusps.  Only  one rational cuspidal curve with
four cusps is known, the quintic with cuspidal configuration [($2_3$); ($2$); ($2$); ($2$)], in other words singularity types $3A_2+A_6$.
It is conjectured that $\kappa \leq 3$ for rational cuspidal curves of degree $d \geq 6$, see \cite{Pion} for more details.

Several classification results for rational cuspidal curves have been obtained and some of them are recalled in the third section. They lead to series of rational cuspidal curves, which have sometimes additional properties, e.g. some of these curves are projectively rigid.

In this note we investigate a stronger property for these series of curves, namely we search among them the free divisors. This is motivated by the fact that the number of known examples of  { \it irreducible}
free divisors seems to be very limited: the Cayley sextic as described in \cite{Sim} and the family of rational cuspidal curves with one cusp
\begin{equation} \label{STfam}
C_d: f_d=y^{d-1}z+x^d+ax^2y^{d-2}+bxy^{d-1}+cy^{d}=0, \  \  \ a \ne 0
\end{equation} 
of type $(d,d-1)$ for $d \geq 5$ described in \cite{ST}. Recently, R. Nanduri has constructed a new family of irreducible free divisors of degree $d \geq 5$, containing the family \eqref{STfam} and consisting of curves having a unique singular point $p$, of multiplicity $m=d-1$, with possibly several branches at $p$, see \cite{N}, especially Remark 2.4. The family \eqref{STfam} as well as the family constructed by Nanduri consist only of rational curves in an obvious way (i.e. the variable $z$ occurs only with exponent one in the corresponding polynomials). For the Cayley sextic the rationality is noted in Remark 3.3 in \cite{Sim}, see also \cite{Se}. 
These examples and those given in the present note suggest the following.
\begin{conj}
\label{conj}
An irrreducible plane curve  of degree $d\geq 2$ which is a free divisor is a rational  curve.
\end{conj}

For free divisors involving several irreducible components, it is easy to construct a free divisor $C$ with at least one irreducible component which is irrational using the recent construction by J. Vall\`es in \cite{JV}, e.g. 
\begin{equation} \label{pen}
C:f=xyz(x^3+y^3+z^3)[(x^3+y^3+z^3)^3-27x^3y^3z^3]=0.
\end{equation}

In the second section we  find new properties of the free divisors in $\PP^2$: a numerical characterization of freeness in Theorem \ref{thmFREE} and Conjecture \ref{conj10}, and  the fact that for a free curve $C:f=0$
its degree $d$ and the total Tjurina number $\tau(C)$ determine the Hilbert function of the graded Milnor algebra $M(f)$, see Theorem \ref{thmHP}. We also give information on the possibile values of the total Tjurina number $\tau(C)$ in terms of the degree $d$, see Theorem \ref{thmWH}. This results also explains why a rational cuspidal curve which is free of degree $d\geq 6$ must have some non weighted homogeneous singularity, i.e. the Jacobian ideal $J_f$  is not of linear type.

In the third section we collect, for the reader's convenience and to fix the notations, some classification results of rational cuspidal plane curves due to various authors.

This classification is used in the forth section to construct new families of irreducible free divisors in $\PP^2$, namely an infinite series in Theorem \ref{thm2ii} of rational curves having two cusps, and the beginning curves in three potentially infinite series in  Example \ref{exprop2i}, Example \ref{exprop3}, Example \ref{exprop4} and Conjecture \ref{conj40}. All these curves have type $(d,m)$ with $m\leq d-2$, hence they are distinct from the curves in \eqref{STfam} or in Nanduri's family in \cite{N}.
We also list all the free rational cuspidal curves of degree 6 in Example \ref{exd=6}.

 The computations of various invariants given in this paper were made using two computer algebra systems, namely CoCoA \cite{Co} and Singular \cite{Sing}, and play a key role especially in the final section.
The corresponding codes are available on request, some of them being available in \cite{St}.

The first author thanks Aldo Conca for some useful discussions.

\section{Free divisors and Milnor algebras} \label{sec2}

Let $f$ be a homogeneous polynomial of degree $d$ in the polynomial ring $S=\C[x,y,z]$ and denote by $f_x,f_y,f_z$ the corresponding partial derivatives.
Let $C$ be the plane curve in $\PP^2$ defined by $f=0$ and assume that $C$ is reduced. We denote by $J_f$ the Jacobian ideal of $f$, i.e. the homogeneous ideal of $S$ spanned by $f_x,f_y,f_z$ and denote by $M(f)=S/J_f$ the corresponding graded ring, called the Jacobian (or Milnor) algebra of $f$. Let $I_f$ denote the saturation of the ideal $J_f$ with respect to the maximal ideal $(x,y,z)$ in $S$. 

 Consider the graded $S-$submodule $AR(f) \subset S^{3}$ of {\it all relations} involving the derivatives of $f$, namely
$$\rho=(a,b,c) \in AR(f)_m$$
if and only if  $af_x+bf_y+cf_z=0$ and $a,b,c$ are in $S_m$. We set $ar(f)_k=\dim AR(f)_k$ and $m(f)_k=\dim M(f)_k$ for any integer $k$.  Then $C$ is a free divisor if $AR(f)$ is a free graded $S$-module of rank two.

Recall  the following basic fact, see \cite{ST}, \cite{Se}, \cite{DS14}.

\begin{prop}
\label{prop5}
Let $C$ be the plane curve in $\PP^2$ defined by $f=0$ and assume that $C$ is reduced and $d=\deg (f)$. Then the following hold.

\noindent (i) The curve $C$ is a free divisor if and only if $I_f=J_f$.

\noindent (ii) The curve $C$ is projectively rigid if and only if the degree $d$ homogeneous components $I_{f,d}$ and $J_{f,d}$ coincide.

In particular, any free divisor is projectively rigid.

\end{prop}

We recall also some definitions, see \cite{DStEdin}.

\begin{definition}
\label{def}

For a plane curve $C:f=0$ of degree $d$ with isolated singularities we introduce three integers, as follows.

\noindent (i) the {\it coincidence threshold} 
$$ct(f)=\max \{q:\dim M(f)_k=\dim M(f_s)_k \text{ for all } k \leq q\},$$
with $f_s$  a homogeneous polynomial in $S$ of degree $d$ such that $C_s:f_s=0$ is a smooth curve in $\PP^2$.

\noindent (ii) the {\it stability threshold} 
$st(f)=\min \{q~~:~~\dim M(f)_k=\tau(C) \text{ for all } k \geq q\},$
where $\tau(C)$ is the total Tjurina number of $C$, that is
$\tau(C)=\sum_{i=1,p} \tau(C,{a_i}).$

\noindent (iii) the {\it minimal degree of a  syzygy} $mdr(f)=\min \{q~~:~~ H^2(K^*(f))_{q+2}\ne 0\}$,
where $K^*(f)$ is the Koszul complex of $f_x,f_y,f_z$ with the natural grading.
\end{definition}
Note that one has  for $j<d-1$ the following equality
\begin{equation} 
\label{ar=er}
AR(f)_j=H^2(K^*(f))_{j+2}.
\end{equation} 
We set $er(f)_j=\dim H^2(K^*(f))_{j+2}$ for any $j$, the number of essential relations among the partial derivatives of $f$.
It is known that one has
\begin{equation} 
\label{REL}
ct(f)=mdr(f)+d-2,
\end{equation} 
It is interesting that the freeness of the plane curve $C$ can be characterized in terms of these invariants. Let $T=3(d-2)$ denote the degree of the socle of the ring $M(f_s)$.

\begin{thm}
\label{thmFREE} For a plane curve $C:f=0$ of degree $d\geq 4$, the following are equivalent.

\noindent (i) $C$ is free;

\noindent (ii)  the 
equality 
$$m(f)_{2d-5-j}+ar(f)_j=\tau(C)$$
holds for any integer $j$ with $-1 \leq j \leq d-2$ and $ar(f)_{d-2} \ne 0$.
In particular
$$st(f)=2d-4-mdr(f)=T-ct(f).$$
\noindent (iii) one has $m(f)_{[\frac{T}{2}]}+m(f)_{T-[\frac{T}{2}]}-m(f_s)_{[\frac{T}{2}]}=\tau(C),$ where $
[\frac{T}{2}]$ denotes the integral part of $\frac{T}{2}$.
\end{thm}

\proof
The fact that the two last formulas in $(i)$ are equivalent follows from \eqref{REL}. Moreover,
if $C$ is free, then $I_f=J_f$, and Proposition 2 in \cite{DBull} implies that
$$\dim M(f) _k= \tau (C)$$
for all $k\geq q$ if and only if $q=T-ct(f)$. The definition of $st(f)$ implies that $st(f)=T-ct(f)$.
Note also that 
$$m(f)_k \geq m(f)_{k+1}\ge \tau(C)$$
 for $k\geq 2d-5$ by Corollary 8 in \cite{CD}. This implies the formula for $st(f)$ given in (ii).

Now we pass to the proof of the theorem, by showing that (i) is equivalent to (ii).
 Consider the rank two vector bundle $T\langle C\rangle=Der(-logC)$ of logarithmic vector fields  along $C$.
By definition, $C$ is free if this bundle splits as a direct sum of two line bundles on $\PP^2$, and this happens exactly when
\begin{equation} 
\label{eq1}
H^1(\PP^2, T\langle C\rangle(k))=I_{f,k+d}/J_{f,k+d}=0
\end{equation} 
for any $k$, see Remark 4.7 in \cite{DS14}.

Using the results in the third section of \cite{DS14}, we know that
$$\chi(T\langle C\rangle(k))=3{k+3 \choose 2}-{d+k+2 \choose 2} +\tau(C),$$
and $h^0((T\langle C\rangle(k))=ar(f)_{k+1}$, $h^2((T\langle C\rangle(k))=ar(f)_{d-5-k}$.

It follows that we have
\begin{equation} 
\label{eq2}
h^1((T\langle C\rangle(k))=ar(f)_{k+1}+ar(f)_{d-5-k}-\tau(C)+{d+k+2 \choose 2}-3{k+3 \choose 2}.
\end{equation} 
 Recall that one has
$$I_{f,j}/J_{f,j}=I_{f,T-j}/J_{f,T-j}$$
for any $j$ by \cite{DS1}, \cite{Se}. Corollary 4.3 in \cite{DPop} shows that in order to check \eqref{eq1}, it is enough to consider only the case
 $$-3 \leq k\leq T/2-d= \frac{d}{2}-3\leq d-4.$$
If we make the substitution $k=d-4-i$, then $0 \leq i \leq d-1$ and the above formula becomes
\begin{equation} 
\label{eq3}
h^1((T\langle C\rangle(k))=ar(f)_{d-3-i}+ar(f)_{i-1}-\tau(C)+{d+i \choose 2}-3{i+1 \choose 2}.
\end{equation} 
Note that both indices $d-3-i$ and $i-1$ are in the interval $[-2,d-2]$.
On the other hand, for $ j \leq d-2$, Theorem 1 in \cite{DBull} and \eqref{ar=er} imply that
\begin{equation} 
\label{eq3.5}
ar(f)_j=m(f)_{d-1+j}-m(f_s)_{d-1+j}.
\end{equation} 
 It follows that
\begin{equation} 
\label{eq4}
ar(f)_{d-3-i}=m(f)_{2d-4-i}-m(f_s)_{2d-4-i}=m(f)_{2d-4-i}-{d+i \choose 2}+3{i+1 \choose 2}.
\end{equation} 
This gives
$$h^1((T\langle C\rangle(k))=m(f)_{2d-4-i}+ar(f)_{i-1}-\tau(C).$$
Assume now that $ar(f)_{d-2} = 0$,  which clearly implies $st(f) \leq d-3$. 
Note that for any plane curve $C:f=0$ of degree $d$ one has $ct(f) \leq st(f)$ if $d$ is even and $ct(f)-1 \leq st(f)$ if $d$ is odd and $st(f)= [T/2].$
It follows that either $ct(f) \leq d-3$, which is a contradiction, or $ct(f)=st(f)+1=d-2$ and 
$d-3=[T/2]$, which is again a contradiction.
This clearly completes the proof of the fact that (i) is equivalent to (ii).

To show that (ii) is equivalent to (iii), use
Corollary 4.3 in \cite{DPop} which shows that it is enough to consider only the case $k={[\frac{T}{2}]}-d$.
Then the equality in (ii) transformed using the formula \eqref{eq3.5} yields the formula in (iii).
\endproof

\begin{rk}
\label{rkproof}
To prove in a quicker (but perhaps more mysterious) way the fact that (i) is equivalent to (ii), one may alternatively use Corollary 4 in \cite{DS1}. This result also implies that for a free divisor $C:f=0$ one has
$$m(f)_{2d-5-j}+er(f)_j=\tau(C)$$
 for any integer $j$.
\end{rk}

In fact, for a free divisor $C:f=0$, the dimensions $m(f)_j$ are completely determined by two invariants, namely the degree $d$ and the total Tjurina number $\tau(C)$. More precisely, one has the following result.

\begin{thm}
\label{thmHP} Let $C:f=0$ be a free divisor of degree $d$ and  total Tjurina number $\tau(C)$, which is not a pencil of lines. Let $d_1$ and $d_2$ with $d_1 \leq d_2$ be the degrees of two homogeneous generators of the free graded $S$-module $AR(f)$. Then the following holds.

\noindent (i) The degrees $d_1$ and $d_2$ are the roots of the equation
$$t^2-(d-1)t+(d-1)^2-\tau(C)=0.$$
In particular $d=d_1+d_2+1$ and $\tau(C)=(d-1)^2-d_1d_2$ and hence the pairs $(d,\tau(C))$ and $(d_1,d_2)$ determine each other.

\noindent (ii) $mdr(f)=d_1$, $ct(f)=d+d_1-2$ and $st(f)=d+d_2-3$.

\noindent (iii)  $ct(f)\leq d+j\leq st(f)$ if and only if   $d_1-2 \leq j\leq d_2-3$, and for such $j$'s one has
$$m(f)_{d+j}=m(f_s)_{d+j}+ {j-d_1+3 \choose 2}.$$
In particular, one has
$$\tau(C)=m(f_s)_{d+d_2-3}+ {d_2-d_1 \choose 2}.$$

\noindent (iv) Let $U=\PP^2 \setminus C$. Then the Euler number $E(U)$ of $U$ is given by
$$E(U)=\tau(C)-\mu(C) +(d_1-1)(d_2-1),$$
where $\mu(C)$ is the total Milnor number of  $C$. In particular, if $C$ is irreducible one has $E(U) \geq 1$ and $d_1>1$.

\end{thm}
Note that for a line arrangement $C$ in $\PP^2$ one has $\tau(C)=\mu(C)$, $b_1(U)=d-1=d_1+d_2$ and $b_2(U)=E(U)+(d-1)-1=d_1d_2$, the claims (i) and (iv) above are a very special case of Terao's results in \cite{Te}. See also \cite{Yo} for a very good survey on free line arrangements.

\proof The definition of the degree $d_1$ and $d_2$ is equivalent to the equality
$$T\langle C\rangle(-1)=\OO(-d_1)\oplus \OO(-d_2).$$
Then the first claim follows from Lemma 4.4 in \cite{DS14} and the second claim follows
from \eqref{REL}, Theorem \ref{thmFREE} and the obvious fact $mdr(f)=d_1$. Note that $d_1>0$ as $C$ is not a pencil of lines. Hence $d_2-2\leq d-4$.
Then equation \eqref{eq3.5}
implies $m(f)_{d+j}=m(f_s)_{d+j}+ar(f)_{j+1}$ for $j \leq d-3$. Since $ j \leq d_2-2$, one has
$$ar(f)_{j+1}= {j-d_1+3 \choose 2},$$
which completes the proof of (iii). To prove (iv), we use the formula
\begin{equation} 
\label{EC}
E(C)=2-(d-1)(d-2) +\mu(C),
\end{equation} 
see for instance \cite{D1}. Then we compute 
$$E(U)=E(\PP^2) -E(C)=3-(2-(d-1)(d-2) +\mu(C))$$
 and this yields the claimed formula using (i). The final claim is clear, since $C$ irreducible is equivalent to $b_1(U)=0$.
\endproof

\begin{rk}
\label{rkRes} The equations $d=d_1+d_2+1$ and $$(d-1)^2-d_1d_2=m(f_s)_{d+d_2-3}+ {d_2-d_1 \choose 2}$$ obtained from Theorem \ref{thmHP} do not impose  restrictions on the integers $1 \leq d_1 \leq d_2$, as they reduce to an identity involving $d_1$ and $d_1$.
\end{rk}

The study of a large number of examples suggests the following conjecture.
\begin{conj}
\label{conj10}
A plane curve $C:f=0$ is free if and only if 
$$ct(f)+st(f)=T.$$
\end{conj}

The following result gives some restrictions on the total Tjurina number of an irreducible  free curve.
\begin{thm}
\label{thmWH} Let $C:f=0$ be an irreducible  free divisor of degree $d>1$. Then $d \geq 5$ and the following hold.

\noindent (i) If $d=5$, then $C$ is rational and all singularities of $C$ are cusps and weighted homogeneous. Moreover in this case $d_1=d_2=2$ and $\tau(C)=12$.

\noindent (ii) If $d\geq 6$, then $C$ is either irrational, or $C$ has at least one singularity which is either not a cusps or it is not weighted homogeneous. Moreover in this case one has
$$\frac{3}{4}(d-1)^2 \leq \tau(C) \leq d^2-4d+7,$$
$d_1\geq 2$ and the integer $\Delta=4\tau(C)-3(d-1)^2$ is a perfect square.
\end{thm}

\proof
 One has
$E(C)=b_0(C)-b_1(C)+b_2(C)=2-b_1(C) \leq 2.$ It follows from \eqref{EC} that
 \begin{equation} 
\label{mu}
 \mu(C) \leq (d-1)(d-2)
\end{equation} 
and the equality holds if and only if $C$ is rational and cuspidal as only then $b_1(C)=0$.
On the other hand, one clearly has $\tau(C) \leq \mu(C)$ with equality if and only if all the singularities of $C$ are weighted homogeneous. Lemma 4.4 in \cite{DS14} imply that
\begin{equation} 
\label{eq10}
4\tau(C)-3(d-1)^2=u^2
\end{equation} 
for  $u=d_2-d_1$. In particular, we get  by putting everything together
\begin{equation} 
\label{eq11}
\frac{3}{4}(d-1) \leq \frac{\tau(C)}{d-1} \leq d-2.
\end{equation} 
This clearly implies $d \geq 5$. If $d=5$, then the two extreme terms coincide to 3, hence we should have equalities everywhere. This proves the claim (i) by using Theorem \ref{thmHP}.

Assume now $d>5$.Then  $d_1 \geq 2$ as shown in Theorem \ref{thmHP} (iv) and  the obvious fact that $\tau(C)=(d-1)^2-d_1d_2$ is maximal (resp. minimal) when the difference $d_2-d_1$ is maximal i.e. when $d_1=2$ (resp. minimal, i.e. when $d_1=[(d-1)/2]$) yields the claimed inequalities.

\endproof

\begin{cor}
\label{corRC}
If $C$ is a free irreducible curve, then $C$ is rational cuspidal if and only if
$$(d_1-1)(d_2-1)=\mu(C)-\tau(C)+1.$$
In particular, a free rational cuspidal curve cannot have only weighted homogeneous singularities unless $d_1=d_2=2$, and hence $d=5$.
\end{cor}

\proof The proof follows from the fact that an irreducible curve $C$ is rational cuspidal if and only if $C$ is homeomorphic to $\PP^1$, and this happens exactly when $E(C)=E(\PP^1)$ which is equivalent to $E(U)=1$.

\endproof

\begin{rk}
\label{rkWH} (i) As shown in Proposition 1.6 in \cite{ST} and mentioned in Remark 4.7 in \cite{DS14}, 
the curve $C:f=0$ has only weighted homogeneous singularities if and only if the Jacobian ideal $J_f$ is of linear type. In particular, a rational cuspidal free divisor of degree $d\geq 6$ cannot have a  Jacobian ideal $J_f$  of linear type by Theorem \ref{thmWH}. This is the case with the family \eqref{STfam}.

(ii) Note that the perfect square $\Delta$ can be zero for $d$ odd, for instance  for the degree $13$ curve $C_1$  described in Proposition \ref{prop4} and Example \ref{exprop4}, which has $\tau(C_1)=108$ and for the curves described in Theorem \ref{thm2ii} which have $d=2k+1$ and $\tau(C)=3k^2$.
In fact, with the notation from Theorem \ref{thmHP}, we have $\Delta=(d_1-d_2)^2$.
For the family \eqref{STfam} one has $\tau(C_d) =d^2-4d+7$ (and this happens for any irreducible free curve with $d_1=2$ as noticed above), hence both inequalities given
in Theorem \ref{thmWH} (ii) for $\tau(C)$ are sharp.

(iii) If $C$ is not irreducible, then one has the following inequalities in analogy to those in Theorem \ref{thmWH}.
Note that $E(C) \leq 2$ is equivalent to $E(U)=E(\PP^2)-E(C) \geq 1$.
Hence in this case we also get $d_1 \geq 2$, and the inequalities in  
Theorem \ref{thmWH} (ii) work unchanged. 

If $E(U)=0$, then it is known that all the irreducible components of $C$ are rational curves, see \cite{GP}, \cite{WV}. In this case we get from $E(C)=3$ the inequality 
$$ \tau (C) \leq \mu(C) = (d-1)(d-2)+1.$$
The equality $ \tau (C) =\mu(C)=(d-1)(d-2)+1$ implies as via Theorem \ref{thmHP} (iv), that $d_1=1$.
Hence the curve $C$ admits a 1-dimensional symmetry group and such curves have been studied in \cite{duW}. In particular, there are free divisors among them, see Proposition 1.3 (2a) in \cite{duW}. Moreover, the corresponding symmetry group is semi-simple if $d \geq 6$ by \cite{duW}, the corresponding degree one syzygy can be diagonalized and hence the corresponding curves are among those studied in the third section of \cite{BC}. The subcase $\tau (C) \leq  (d-1)(d-2)$ leads again to $d_1 \geq 2$, and the inequalities in  
Theorem \ref{thmWH} (ii) work unchanged. 

Finally, if $E(U)<0$, then any irreducible component $C_i$ is a rational cuspidal curve and there is a unique point $p$ such that $C_i \cap C_j =\{p\}$ for any $i \ne j$, i.e. the curve $C$ is very special. We do not know which free divisors occur in this way, except the union of lines passing through one point.

In conclusion, with the exception of the special cases listed here, one has
$ \tau (C) \leq d^2-4d+7$ even for a reducible free divisor $C$.

\end{rk}

\begin{ex}\label{exd=7,8}
Let $C$ be a free divisor of degree $d=6$ which is rational and cuspidal. It follows from Theorem \ref{thmWH} that $C$ has some non weighted homogeneous singularity and moreover
$\tau(C)=19.$
Let $C$ be a free divisor of degree $d=7$ which is rational and cuspidal. It follows from Theorem \ref{thmWH} that
$\tau(C) \in \{27,28\}$ . However, all the examples of such curves coming from  \cite{FZ},  \cite{SaTo},  \cite{ST}  have $\tau(C) =28.$ On the other hand, there are free line arrangements of degree 7 having $\tau=27$, for instance the line arrangement whose defining polynomial is
$$f=(x+z)(x-z)(y+z)(y-z)(x+y)(x-y)z.$$
Similarly, let $C$ be a free divisor of degree $d=8$ which is rational and cuspidal. It follows from Theorem \ref{thmWH} that
$\tau(C) \in \{37,39\},$  but all the examples of such curves coming from  \cite{FZ},  \cite{SaTo},  \cite{ST}  have $\tau(C) =39.$  But there are free line arrangements of degree 8 having $\tau=37$, for instance that defined by the polynomial
$$f=(x+z)(x-z)(y+z)(y-z)(x+y)(x-y)yz.$$
However, for $d=9$ it follows from Theorem \ref{thmWH} that
$\tau(C) \in \{48,49, 52\},$  and   the divisor $C_9$ in Remark \ref{rkeq2i} has $\tau(C_9)=52$,  while the curve $C_9$  described in Theorem \ref{thm2ii} has $\tau(C_1)=48$. Moreover, the free line arrangement given by
$$f=(x+z)(x-z)(y+z)(y-z)(x+y)(x-y)xyz=0$$
has $\tau(C)=49$.
\end{ex}

\section{Some classification results for  rational cuspidal curves} 

In this section we recall some of the classification results for rational cuspidal curves in $\PP^2$.
First we consider the case of rational cuspidal curves of type $(d,d-2)$, following the work by Flenner-Zaidenberg \cite{FZ} and Sakai-Tono \cite{SaTo}.  Note that a rational cuspidal curve of type $(d,d-2)$ may have as its "largest" cusp a cusp with the unique Puiseux characteristic pair $(d-1,d-2)$, i.e. topologically equivalent to the cusp $u^{d-1}+v^{d-2}=0$, see Examples \ref{exprop2i} and \ref{exprop3} below. Such a cups is called in the sequel a cusp of type $(d-1,d-2)$. However, this largest cusp can also be a cusp with the unique Puiseux characteristic pair $(d,d-2)$ for $d$ odd, i.e. topologically equivalent to the cusp $u^{d}+v^{d-2}=0$, called a cusp of type $(d,d-2)$, see Theorem \ref{thm2ii} for  an infinite series of examples.

In increasing number of cusp order, we have the following results.
The case of one cusp is quite clear-cut, and it is described in \cite{SaTo}. 
\begin{prop}
\label{prop1}
Let $C$ be a rational cuspidal curve of type $(d,d-2)$ having a unique cusp. Then $d$ is even and up to projective equivalence the equation of $C$ can be written as
$$f=(y^kz+\sum_{i=1,k+1}a_ix^i y^{k+1-i})^2 - xy^{2k+1}=0,$$
where $a_{k+1}\ne 0$ and $d=2k+2 \geq 4$.
\end{prop}

The case of two  cusps is more involved, and is described in \cite{SaTo}. 
\begin{prop}
\label{prop2}
Let $C_d:f_d=0$ be a rational cuspidal curve of type $(d,d-2)$ having two cusp, let's say $q_1$ of multiplicity $d-2$ and $q_2$ of multiplicity $\leq d-2$. Then one of the following three cases arises.

\noindent (i) The germ $(C_d,q_2)$ is a singularity of type $A_{2d-4}$, hence of multiplicity $2$ and Milnor (or Tjurina) number $2d-4$. Then for each $d \geq 4$, $C_d$ is unique up to projective equivalence.

\noindent (ii) The germ $(C_d,q_2)$ is a singularity of type $A_{d-1}$, with $d$ odd. Up to projective equivalence the equation of $C_d$ can be written as
$$C:f=(y^{k-1}z+\sum_{i=2,k}a_ix^i y^{k-i})^2y - x^{2k+1}=0,$$
where $d=2k+1 \geq 5$.

\noindent (iii) The germ $(C_d,q_2)$ is a singularity of type $A_{2j}$, with $d$ even and $1 \leq j \leq (d-2)/2$.  Up to projective equivalence the equation of $C_d$ can be written as
$$C_d:f_d=(y^{k+j}z+\sum_{i=2,k+j+1}a_ix^i y^{k+j+1-i})^2 - x^{2j+1}y^{2k+1}=0,$$
where $a_{k+j+1} \ne 0$, $d=2k+2j+2 \geq 6$,  $k \geq 0$, $j \geq 1$.
\end{prop}

\begin{rk}
\label{rkeq2i}
To obtain the equations $f_d=0$ for the curves $C_d$ in Proposition \ref{prop2}  $(i)$, one can proceed as follows. Start with the equation for the cuspidal cubic $C_3: f_3(x,y,z)=yz^2-x^2z+x^3$. For $d \geq 4$, to get the polynomial $f_d$ from the previous polynomial $f_{d-1}$,
supposed already computed, we proceed as follows.
First we perform the substitution in $f_{d-1}$ given by
$$x \mapsto x^2,  \\  y \mapsto xy, \\  z \mapsto yz+a_{d-1}x^2,$$
with $a_{d-1}$ the coefficient of $x^{d-1}$ in $f_{d-1}$, then we divide by $x^{d-3}y$ and denote the resulting polynomial by $f_d$. In other words
$$f_d(x,y,z)=f_{d-1}(x^2,xy, yz+a_{d-1}x^2)x^{-d+3}y^{-1}.$$
We get in this way
$$C_4: f_4(x,y,z)=(yz+x^2)^2-x^3z=0,$$
$$C_5: f_5(x,y,z)=y(yz+x^2)^2+2x^3(yz+x^2)-x^4z=0,$$
$$C_6:(yz + 2x^2)^2y^2 + 2(yz + 2x^2)yx^3 + (2yz + 5x^2)x^4 - x^5z = 0,$$
These equations occur already in  \cite{SaTo}. 
The next few curves in this family are  the curves
$$C_7:f_7=14x^7 + 14x^6y + 20x^5y^2 + 25x^4y^3 - x^6z + 2x^5yz + 2x^4y^2z + 4x^3y^3z + 10x^2y^4z $$
$$+ y^5z^2=0,$$
$$C_8: f_8= 42x^8 + 48x^7y + 81x^6y^2 + 140x^5y^3 + 196x^4y^4 - x^7z + 2x^6yz + 2x^5y^2z $$
$$+ 4x^4y^3z + 10x^3y^4z + 
28x^2y^5z + y^6z^2=0,$$
$$C_9: f_9=132x^9 + 165x^8y + 308x^7y^2 + 616x^6y^3 + 1176x^5y^4 + 1764x^4y^5 - x^8z + 2x^7yz $$
$$+ 2x^6y^2z + 4x^5y^3z + 
10x^4y^4z + 28x^3y^5z + 84x^2y^6z + y^7z^2=0,$$
 and finally
$$C_{10} :f_{10}=429x^{10} + 572x^9y + 1144x^8y^2 + 2496x^7y^3 + 5460x^6y^4 + 11088x^5y^5 + 17424x^4y^6 $$
$$- x^9z + 2x^8yz + 
2x^7y^2z + 4x^6y^3z + 10x^5y^4z + 28x^4y^5z + 84x^3y^6z + 264x^2y^7z + y^8z^2=0.$$
The curve $C_d:f_d=0$ for $5 \leq d \leq 15$ has $d_1=2$, the corresponding relation being
$$A_df_{d,x}+B_df_{d,y}+C_df_{d,z}=0,$$
with the coefficients $A_d= (d-2)x^2+4(d-3)xy$, $B_d=2(d-1)xy-4(2d-3)y^2$, $C_d=2d(2d-7)a_{d-1}x^2-(d-1)(d-2)xz+2(d-2)(2d-3)yz$, and the coefficient $a_{d-1}$ as introduced above. It follows from Theorem \ref{thmHP} that $d_2=d-3$ and $\tau(C_d)=d^2-4d+7$.

\end{rk}

Now consider the case of three cusps, considered in \cite{FZ}. One has the following result.

\begin{prop}
\label{prop3}
Let $C$ be a rational cuspidal curve of type $(d,d-2)$ having three cusps. Then there exists a unique pair of integers $a,b$, $a\geq b \geq 1$ with $a+b=d-2$ such that up to projective equivalence the equation of $C$ can be written in affine coordinates $(x,y)$ as
$$f(x,y)= \frac{x^{2a+1}y^{2b+1}-((x-y)^{d-2}-xyg(x,y))^2}{(x-y)^{d-2}},$$
where $d \geq 4$,  $g(x,y)=y^{d-3}h(x/y)$ and 
$$h(t)= \sum_{k=0,d-3}\frac{a_k}{k!}(t-1)^k,$$
with $a_0=1$, $a_1=a-\frac{1}{2}$ and $a_k=a_1(a_1-1) \cdots (a_1-k+1)$ for $k>1$.
\end{prop}

\begin{rk}
\label{rkrigid}
The above results imply that the  cuspidal rational curves in Proposition \ref{prop2} $(i)$ and Proposition \ref{prop3} are projectively rigid, a fact explicitly discussed in \cite{FZ}.
\end{rk}

Another classification problem is to find all the Puiseux pairs $(a,b)$ such that there is a rational unicuspidal (i.e. $\kappa=1$) curve $C$ in $\PP^2$, whose cusp has a unique Puiseux pair $(a,b)$. The complete result is described in Theorem 1.1 in \cite{FLMN}, but here we just mention the cases  $(b)$ and $(d)$ of this Theorem. See also \cite{BL}, section 2.3 and   \cite{Ka}, Corollary 11.4.

\begin{prop}
\label{prop4}

(i) The curve $C_d:f_d=(zy-x^2)^k-xy^{2k-1}=0$ is a rational unicuspidal plane curve  of degree $d=2k$ for $k\geq 2$, with a unique Puiseux pair $(k,4k-1)$.

(ii) Let $a_i$ be the Fibonacci numbers with $a_0=0$, $a_1=1$, $a_{j+2}=a_{j+1}+a_{j}$. Then there is a rational unicuspidal plane curve $C_k$ of degree $a_{2k+5}$ for $k\geq 0$, with a unique Puiseux pair $(a_{2k+3},a_{2k+7})$ and whose defining affine equation is obtained as follows.
Set $P_{-1}=y-x^2$, $Q_{-1}=y$, $P_0=(y-x^2)^2-2xy^2(y-x^2)+y^5$, $Q_0=y-x^2$,
$G=xy-x^3-y^3$. Then define for any $k >0$ recursively
$Q_k=P_{k-1}$ and $P_k=(G^{a_{2k+3}}+Q_k^3)/Q_{k-1}$. Then $P_k=0$ is the defining affine equation for the curve $C_k$.
\end{prop}

\section{New examples of irreducible free divisors} \label{sec3}

Motivated by this Proposition \ref{prop5} and by Remark \ref{rkrigid}, as well as the discussion in the Introduction, we now search for free divisors among the cuspidal rational curves listed in the previous section.
We start with the rigid curves, and recall that a rational quartic is never a free divisor, see \cite{ST}.

\begin{ex}\label{exprop2i}
The curves $C_d$, with $5 \leq d \leq 15$ from Proposition \ref{prop2} $(i)$, whose equations are given in 
Remark \ref{rkeq2i} are free divisors. Indeed, the condition $I_{f,d}=J_{f,d}$ is in fact equivalent to the condition $\dim I_{f,d}=\dim J_{f,d}$ as we have by definition $J_{f,d}\subset I_{f,d}$.
And the equality of dimensions is easily checked using Singular, by computing the Hilbert functions associated to the graded rings $S/I_f$ and $M(f)$.
Each of these curves $C_d$ has two cusps, one of type $(d-1,d-2)$ and the other of type $(2,2d-3)$, i.e. a singularity $A_{2d-2}$. Moreover $d_1=2$ for all of them as shown in Remark \ref{rkeq2i}, and hence 
$\tau(C_d)=d^2-4d+7$ and $\mu(C_d)= (d-2)(d-3)+2d-2=d^2-3d+4$ which is compatible with the formula in Corollary \ref{corRC}.
\end{ex}

\begin{ex}\label{exd=6}
 The cuspidal rational curves of degree 6 have been classified by Fenske, \cite{F},
who obtained 11 main cases. 

The 1-cuspidal curves in this classification fall into 3 classes: $C_1$ (with subcases (a), (b), (c), (d)), $C_2$ (with subcases (a), (b)), $C_3$ (with subcases (a), (b), (c)).
Among these families, only the type $C_3$ (a) is free. It follows that the curve obtained for $d=6$ in the family \eqref{STfam} belongs to this class.

The 2-cuspidal curves fall into 6 classes: $C_4$ (with subcases (a), (b)), $C_5$, $C_6$, $C_7$, $C_8$ and $C_9$. Among these families, the types $C_4$ (b), $C_5$, $C_6$ and $C_9$ are free.

The two classes of 3-cuspidal curves are special cases of the curves discussed  in Proposition \ref{prop3} and Example \ref{exprop3}

All the free divisors above satisfy $ct(f)=st(f)=6$ and $\tau(C)=19$. To decide which divisors are free we use Singular computations, and for families involving parameters we use  Lemma 1.1 in \cite{ST}.
\end{ex}

\begin{ex}\label{exprop3}
The curves $C_d$  from Proposition \ref{prop3} for $5 \leq d \leq 10$ are free divisors.
This is proved by following the same approach as in the previous example. Moreover, in this case we have $d_1=2$ and hence again $\tau(C_d)=d^2-4d+7$. The description of singularities in 
Proposition \ref{prop3} implies that $\mu(C)= d^2-3d+2$ for any $(a,b)$.

\end{ex}

\begin{ex}\label{exprop4}

(i) The curves $C_d$  from Proposition \ref{prop4}  (i) for $6 \leq d=2k \leq 20$ are free divisors with $d_1=2$ and hence $\tau(C_d)=d^2-4d+7$.

(ii) The curves $C_k$  from Proposition \ref{prop4}  (ii) for $0 \leq k \leq 3$ are free divisors and have the following invariants.

\noindent (0) The curve $C_0$ has degree $d=5$, $\tau(C)=12$, $\mu(C)=12$ and $d_1=2$. 

\noindent (1) The curve $C_1$ has degree $d=13$, $\tau(C)=108$, $\mu(C)=132$ and $d_1=6$. 

\noindent (2) The curve $C_2$ has degree $d=34$, $\tau(C)=823$, $\mu(C)= 1056$ and $d_1=14$. 

\noindent (3) The curve $C_3$ has degree $d=89$, $\tau(C)=5889$, $\mu(C)=7656$ and $d_1=35$. 

This is proved by following the same approach as in the previous example. 
\end{ex}

It is natural to suggest the following conjecture.
\begin{conj}
\label{conj40}
Any of the  rational cuspidal curves described in Proposition \ref{prop2} $(i)$, Proposition \ref{prop3}, and  Proposition \ref{prop4} is a free divisor if its degree $d$ is at least $5$.
\end{conj}

Next we present an infinite series of irreducible free divisors obtained from the classification of 
rational cuspidal curves. This family of curves is a special case of the family described in Proposition \ref{prop2} $(ii)$.

\begin{thm}
\label{thm2ii}
The rational cuspidal curve
$$C_{2k+1}: f_{2k+1}=(y^{k-1}z+x^k)^2y-x^{2k+1}=0$$
of type $(2k+1,2k-1)$ has two cusps of type $(2k+1,2k-1)$ and respectively $(2k+1,2)$ and is a free divisor for any $k \geq 2$. The corresponding Jacobian ideal $J_{f_{2k+1}}$ is of linear type if and only if $k=2$.
Moreover  $\tau( C_{2k+1})=3k^2$,   $\mu( C_{2k+1})=2k(2k-1)$ and $d_1=d_2=k$.

\end{thm}

\proof  To prove that we have free divisors for any $k$ we proceed as follows.
We look at the syzygies among the partial derivatives $f_x,f_y,f_z$ and find that we have two such syzygies in degree $k$, namely:
$$(r_1 ): \  \  \ a_xf_x+a_yf_y+a_zf_z=0,$$
where
$a_x=2x^k+2y^{k-1}z$, 
       $a_y=(4k+2) x^k -4k x^{k-1}y-(8k^2-2) y^{k-1}z$,\\
$ a_z=4k(k-1) x^{k-1}z + (8k^3-4k^2-2k+1) y^{k-2}z^2$ and
$$(r_2 ): \  \  \ b_xf_x+b_yf_y+b_zf_z=0,$$
where
$b_x=0$,   
       $b_y=-2y^k$ and  
       $b_z=x^k+(2k-1)y^{k-1}z$. It is clear that $(r_1)$ and $(r_2)$ are linearly independent as $b_x=0$ and $a_x \ne 0$. Then we apply Lemma 1.1  and Proposition 1.8 in \cite{ST}, exactly as in the proof of Proposition 2.2 in \cite{ST}.

The Milnor number $\mu( C_{2k+1})$ is computed using Corollary \ref{corRC} and this implies that the largest cusp of $C_{2k+1}$ has type $(2k+1,2k-1)$. Indeed, the other cusp is described in 
Proposition \ref{prop2} $(ii)$ and we know that it has type $(2k+1,2)$, i.e. it is an $A_{2k}$-singularity with Milnor number $2k$.
\endproof

\end{document}